\documentclass[english,american]{amsart}
\usepackage[T1]{fontenc}
\usepackage[latin1]{inputenc}
\usepackage{graphicx}
\usepackage{amssymb}

\makeatletter
 \theoremstyle{plain}
 \theoremstyle{plain}    
 \newtheorem{thm}{Theorem} 
 \theoremstyle{plain}    
 \newtheorem{lem}{Lemma} 
 \theoremstyle{plain}    
 \newtheorem{cor}{Corollary} 
 \theoremstyle{plain}    
 \newtheorem{prop}{Proposition} 
 \newcommand{\lyxrightaddress}[1]{
   \par {\raggedleft \begin{tabular}{l}\ignorespaces
   #1
   \end{tabular}
   \vspace{1.4em}
   \par}
 }

\usepackage[mathscr]{eucal}
\usepackage{amssymb}
\makeatother

\usepackage{babel}
\makeatother
\begin{document}
\newcommand{\Abb}{\mathbb{A}} \newcommand{\NN}{\mathbb{N}} \newcommand{\RR}{\mathbb{R}} \newcommand{\ZZ}{\mathbb{Z}} \newcommand{\CC}{\mathbb{C}} \newcommand{\QQ}{\mathbb{Q}} \newcommand{\TT}{\mathbb{T}} \newcommand{\PP}{\mathbb{P}} \newcommand{\EE}{\mathbb{E}} \newcommand{\HH}{\mathbb{H}} \newcommand{\sA}{{\mathscr{A}}} \newcommand{\sB}{{\mathscr{B}}} \newcommand{\sE}{{\mathscr{E}}} \newcommand{\sF}{{\mathscr{F}}} \newcommand{\sG}{{\mathscr{G}}} \newcommand{\sH}{{\mathscr{H}}} \newcommand{\sJ}{{\mathscr{J}}} \newcommand{\sK}{{\mathscr{K}}} \newcommand{\sP}{{\mathscr{P}}} \newcommand{\sL}{{\mathscr{L}}} \newcommand{\sT}{{\mathscr{T}}} \newcommand{\sY}{{\mathscr{Y}}} \newcommand{\sZ}{{\mathscr{Z}}} \newcommand{\fN}{\mathfrak{N}} \newcommand{\fA}{\mathfrak{A}} \newcommand{\fS}{\mathfrak{S}} \newcommand{\fZ}{\mathfrak{Z}} \newcommand{\fp}{\mathfrak{p}} \newcommand{\fu}{\mathfrak{u}} \newcommand{\fs}{\mathfrak{s}} \newcommand{\id}{\mathrm{Id}} \newcommand{\ds}{\displaystyle} \newcommand{\limSup}[1][n]{\overline{\lim_{#1\to\infty}}} \newcommand{\iso}{\textrm{Iso}(\TT)}
\newcommand{\esp}[1]{\EE\left[ #1 \right] }

\newcommand{\pr}[1]{\PP\left[ #1 \right] }

\newcommand{\Rp}{\RR_{+}^{*}}
 
\newcommand{\RxR}{\Rp\times\RR}

\newcommand{\affR}{\mathrm{Aff}(\RR)}

\newcommand{\affT}{\mathrm{Aff}(\TT)}

\newcommand{\affO}[1]{\mathrm{Aff}(#1 )}
 
\newcommand{\affq}{\mathrm{Aff}(\QQ_{p})}

\newcommand{\affQ}{\mathrm{Aff}(\QQ)}

\newcommand{\hor}{\mathrm{Hor}(\TT)}

\newcommand{\bT}{\partial^{*}\TT}
 
\newcommand{\bG}{\partial^{*}\Gamma}

\newcommand{\horg}{\mathrm{Hor}(\Gamma)}
 
\newcommand{\horO}[1]{\mathrm{Hor}(#1 )}

\newcommand{\toinf}{\rightarrow\infty}
 
\newcommand{\inv}{^{-1}}
 
\newcommand{\Ind}[1]{1_{#1 }}
 
\newcommand{\nor}[1]{\left| #1 \right| }

\newcommand{\norp}[1]{\left| #1 \right| _{p}}
 
\newcommand{\sumzi}[1]{\sum_{#1 =0}^{\infty}}
 
\newcommand{\dnor}[1]{\left\Vert #1 \right\Vert }

\newcommand{\tm}{\frac{3}{2}}

\newcommand{\stard}{\stackrel{\cdot}{*}}

\newcommand{\supp}{\mathrm{supp}}

\newcommand{\cf}{\wedge}

\newcommand{\e}{\mathrm{e}}

\newcommand{\ver}{\rightarrow}

\newcommand{\Qp}{\QQ_{p}}

\newcommand{\nup}[1]{\nu_{p}\left( #1 \right) }

\newcommand{\nad}[1]{\left\langle #1 \right\rangle }

\newcommand{\sumP}{\sum_{p\in\sP}}

\newcommand{\card}[1]{\textrm{card}\left\{  #1 \right\}  }

\newcommand{\ZP}{\ZZ(P)}

\newcommand{\affP}{\affO{P}}

\newcommand{\gP}{\left( P\right) }

\begin{flushright}\today{}\end{flushright}
\medskip{}

\title{Poisson boundary for finitely generated groups of rational affinities.}

\maketitle
\medskip{}
\begin{center}Sara Brofferio\end{center}
\bigskip{}

\begin{abstract}
The group of affine transformations with rational coefficients, $\affQ$,
acts naturally on the real line $\RR$, but also on the $p$-adic
fields $\Qp$. The aim of this note is to show that all these actions
are necessary and sufficient to represent bounded $\mu$-harmonic
functions for a probability measure $\mu$ on $\affQ$ that is supported
by a finitely generated sub-group, that is to describe the Poisson
boundary. \foreignlanguage{english}{}
\end{abstract}

\section*{Introduction}

On a generic discrete group $G$ there is a natural probabilistic
analogue of classical harmonic functions on differential structures:
given a probability measure $\mu$ on $G$, a bounded measurable function
$f$ is said to be \emph{$\mu$-harmonic} when it satisfies the mean
value property\[
f(x)=\int_{G}f(xg)\mu(dg),\]
that is, the value in a point is always the mean of the values in
its neighborhood according to the measure $\mu$. As in the classical
case a harmonic function is ruled by its behavior on the boundary
of the definition set; the question is to provide a boundary of the
group large enough to represent all harmonic functions. The \emph{Poisson
boundary} is the measurable $G$-space, $B$, endowed with a $\mu$-invariant
probability measure, $\nu$, that enables to represent every bounded
harmonic function, $f$, by the formula \begin{equation}
f(g)=\int_{B}\psi(g\cdot x)\nu(dx)\qquad\forall g\in G\label{eq-front-funz-Intro}\end{equation}
where $\psi$ is a bounded function on $B$.

We are interested here in finding the Poisson boundary for the group
of affine transformations with rational coefficients\[
\affO{\QQ}=\left\{ (a,b):x\mapsto ax+b\,|\, a\in\QQ^{*},\, b\in\QQ\right\} .\]
 This discrete group has a natural action on the real line $\RR$
and a natural injection into the group of real affine transformations.
One can obtain interesting results concerning the behavior of the
probabilistic objects on $\affO{\QQ}$ using the powerful theory developed
on Lie groups, whenever no continuity hypothesis on the measure is
required (for instance \cite{Ke73}, \cite{BP92}, \cite{BBE} or
\cite{Br03b}). However the real line is not the only natural completion
of the rational numbers: others are provided by the the $p$-adic
rationals, $\QQ_{p}$, for any prime number $p$ and one has natural
injection of rational affinities in every group of $p$-adic affine
transformations. Algebraic groups over local fields are popular in
different mathematic fields and have attracted a growing interest
also in probability (cf. for instance \cite{CKW},\cite{AK94} or
\cite{Eva93}). 

The real and the $p$-adic fields formally behave in a very similar
way and it is useful to unify the notations associating the Euclidean
setting with the {}``prime number'' $p=\infty$; thus $\QQ_{\infty}$
is $\RR$, the Euclidean norm is $\nor{\,\cdot\,}_{\infty}$ and so
on. Each of these fields can provide a useful boundary for the study
of harmonic function. In fact, under first moment conditions, whenever
the action of $\affQ$ is contracting in mean on a field $\Qp$, that
is the drift \[
\phi_{p}=\int_{\affQ}\ln\norp{a}d\mu(a,b)\]
 is negative, there is a unique $\mu$-invariant probability measure
$\nu_{p}$ on $\QQ_{p}$ and for any bounded function $\psi$ on $\QQ_{p}$
one can construct a harmonic function on the group using the formula
(\ref{eq-front-funz-Intro}). 

However none of these actions alone is sufficient to represent \emph{all}
bounded harmonic functions. The aim of this note is to show that,
when the measure $\mu$ is supported by a finitely generated sub-group,
the construction of the Poisson boundary requires to consider all
these actions simultaneously. We prove the following

\begin{thm}
\label{thm-poisson-boundary}If $\mu$ is supported by a finitely
generated subgroup of $\affQ$ and has a finite first moment, then
there exists a unique $\mu$-invariant probability measure $\nu$
on \[
B=\prod_{p:\phi_{p}<0}\Qp\]
and $(B,\nu)$ is the Poisson Boundary. Thus there exist non-trivial
harmonic functions if and only if $\phi_{p}\neq0$ for some $p$.
\end{thm}
This is a generalization of the known result on the Baumslag-Solitar
group $BS(1,p)$, whose Poisson boundary was shown by Kaimanovich
and Vershik \cite{KV83} to be either $\RR$ or $\QQ_{p}$ according
to the sign of the drift . We shall adapt their techniques based on
the estimation of the entropy of the tail conditional probability
to our case where there is no more dichotomy between the real and
$p$-adic norms, since one has more degrees of freedom. 

A similar result holds also for some measures whose support generates
the whole group $\affQ$, and not only a finitely generated sub-group.
In this global situation the Poisson Boundary can be proved to be
a sub-space of the Adele ring. We have chosen here to focus on the
simplest situation and to postpone the study of more general cases
to future work in order to allow a more graphical interpretation,
that is not available in infinite dimensional situations. Here, we
can also avoid some heavy technical details, that might obscure the
basic ideas and reduce accessibility.

This paper is structured as follows. In Section 1, we first introduce
the concept of Poisson boundary and $\mu$-boundaries and their relationship
with random walks on groups. We then define the algebraic structures
we are dealing with: the finitely generated sub-groups of $\affQ$
and their geometrical boundary. In Section 2 we prove that all the
fields $\Qp$ with a contractive action contribute to the representation
of harmonic functions. In Section 3 we show that these actions are
sufficient to describe all bounded harmonic functions, that is we
determine the Poisson boundary and prove Theorem 1.

\section{Preliminaries}

\subsection{Poisson boundary}

Given a probability measure $\mu$ on a discrete group $G$ the \emph{Poisson
boundary} is a $G$-space, $B$, endowed with a $\mu$-invariant probability
measure $\nu$ such that every bounded harmonic function $f$ on $G$
can be written in the form \[
f(g)=\int_{B}\psi(g\cdot x)\nu(dx)\qquad\forall g\in G\]
for a (unique) bounded measurable function $\psi$ on $B$ (cf. Furstenberg
\cite{Fur63}).

This space is unique from a measure theoretic point of view and can
be identified with the exit boundary of the associated random walk
on the group. The \emph{random walk} of law $\mu$ is the Markov chain
$\left\{ R_{n}\right\} _{n\in\NN}$ on $G$ with transition kernel
$Pf(x)=\int f(xg)\mu(dg)$, which can also be seen as the process
obtained by iterated multiplication of a sequence $\left\{ g_{n}\right\} _{\in\NN}$
of independent and identically distributed random elements of $G$
with law $\mu$. Then the random walk starting at the identity is
defined by the equations\[
R_{0}=e\quad\,\textrm{and }\quad R_{n+1}=R_{n}g_{n+1}\quad\forall n\in\NN.\]
The law of this process on the space $G^{\NN}$ is denoted by $\PP$
and $\EE$ is the related mean. One says that two paths on $G$ have
the same orbit \emph{}up to a time shift if \[
\left(x_{n}\right)_{n}\sim\left(x'_{n}\right)_{n}\Longleftrightarrow\exists k,h\in\NN:x_{n+k}=x'_{n+h}\forall n\in\NN.\]
The Poisson boundary \emph{}coincides with the quotient of the $(G^{\NN},\PP)$
by this equivalence relation and from this point of view it can be
interpreted as the space that best describes the asymptotic behavior
of the trajectories of the random walk.

In typical situation the group $G$ acts continuously on a topological
space $X$ and $R_{n}\cdot x$ converges almost surely to a random
element whose law $\eta$ is indipendent from $x\in X$. Then the
measure space $(X,\eta)$ is a quotient of the Poisson boundary, that
is a \emph{$\mu$-boundary} in the sense of Furstenberg. To give a
geometrical interpretation of the Poisson boundary one needs to provide
a geometrical space where the random walk converges but that separates
completely the tails of the trajectories.

We now want to describe a geometrical interpretation of finitely generated
sub-groups, $G$, of $\affQ$ and a associated $G$-space that will
be proved to have these two fundamental properties.

\subsection{Finitely generated sub-groups of $\affQ$}

First observe that, as we deal with a finitely generated sub-group,
there exists a finite set $P$ of prime numbers such that any element
$(a,b)\in G$ may be written in the form \[
a=\pm\prod_{p\in P}p^{k_{p}}\,\textrm{ and }b=h_{0}\prod_{p\in P}p^{h_{p}}\qquad\textrm{with }h_{0},h_{p},k_{p}\in\ZZ.\]
The finite set $P$ contains all prime numbers that are significant
in the study of the random walk on $G$. If we denote by $\gP$ the
sub-group of the multiplicative group $\QQ^{*}$ generated by $P$
and $-P$ and denote by $\ZP$ the additive sub-group of the group
$\QQ$ consisting of the numbers of the form $h\prod_{p\in P}p^{h_{p}}$,
we can restrict our study, without loss of generality, to the case
where $G$ is the sub-group \[
\affP=\left\{ (a,b)\in\affQ|a\in\gP,b\in\ZP\right\} =\gP\ltimes\ZP.\]
 Observe that if $P$ is reduced to a singleton, then $\affO{\left\{ p\right\} }$
is an index-2-extention of the Baumslag-Solitar group $BS(1,p)$,
since in our case the linear coefficient $a$ can be positive or negative.

Since $\affP$ is the semi-direct product of $\gP$ and $\ZP$, the
projection $(a,b)\mapsto a$ is a homomorphism on $\gP$. The structure
of this group is quite simple; for any $p\in P$ the \emph{$p$-adic
valuation} $\nu_{p}$, that is the map \[
\nu_{p}:a=\pm\prod_{q\in P}p^{k_{q}}\mapsto k_{p},\]
is a homomorphism onto $\ZZ$ and $\gP$ is a direct sum of \emph{$\card{P}$}
copies of $\ZZ$ and of a copy of $\ZZ_{2}$ to take care of the sign
\[
\gP\cong\ZZ_{2}\bigoplus_{p\in P}\ZZ.\]
 This decomposition gives immediately a discrete geometrical structure
for the dilation-contraction component of our group. 

Let us now focus our attention on the \emph{translation sub-group}
$\ZP$. It can be naturally embedded in the real line, but, since
it is not finitely generated, this embedding is not discrete with
respect with to the Euclidean metric ($\ZP$ is in fact dense in $\RR$).
Thus to provide a structure that completely separates the trajectories
in the group, we have to consider also the embedding of $\ZP$ into
the $p$-adic fields $\QQ_{p}$, that is the completion of $\QQ$
according to the \emph{$p$-adic norm} \[
\left|q\right|_{p}=p^{-v_{p}(q)},\]
where the \emph{$p$-adic valuation} of an integer $r$ is $v_{p}(r)=\max\left\{ k\in\NN|p^{-k}r\in\NN\right\} $
and $v_{p}(r/s)=v_{p}(r)-v_{p}(s).$ 

To unify the notation in $p$-adic and real settings, we are going
to associate, as often occurs in similar context, the real objects
with the symbol $p=\infty$; thus $\QQ_{\infty}$ is $\RR$, the Euclidean
norm is $\nor{\,\cdot\,}_{\infty}$ and so on. We also denote by $\overline{P}$
the union of $P$ and $\infty$. Then one can easily check the following
lemma, that is going to be frequently used in the following sections. 

\begin{lem}
\label{lem-appro-raggioI}The diagonal embedding of $\ZP$ in $\prod_{p\in\overline{P}}\QQ_{p}$is
discrete and for every $(z_{p})_{p\in\overline{P}}\in\prod_{p\in\overline{P}}\QQ_{p}$
the set \[
\left\{ b\in\ZP|\textrm{ }\norp{z_{p}-b}\leq1\forall p\in\overline{P}\right\} \]
 has always 2 or 3 elements.
\end{lem}
The first consequence is that we can see $\affP$ as a discrete sub-set
of $\gP\times\prod_{p\in\overline{P}}\QQ_{p}$ using the injection\begin{eqnarray*}
\affP & \longrightarrow & \gP\times\prod_{p\in\overline{P}}\QQ_{p}\\
(a,b) & \mapsto & (a,(b)_{p\in\overline{P}}).\end{eqnarray*}
Obviously, this map is also an homomorphism when one considers $\gP\times\prod_{p\in\overline{P}}\QQ_{p}$
as a semi-direct product with respect of the action induced by the
usual multiplication of $\gP$ on $\Qp$. Furthermore, if $p$ is
a {}``true'' prime ($p\neq\infty$), it is the $p$-adic valuation
of an element $a\in\gP$ that determines whether this action is dilating
or contracting; in fact\[
\norp{ax}=\norp{a}\norp{x}=p^{-v_{p}(a)}\norp{x}.\]
 Thus the decomposition of $\gP$ as direct sum of different copies
of $\ZZ$ corresponds to the contracting directions of the actions
on the $\QQ_{p}$ for the each different $p\in P$. On the other hand,
the Euclidean norm is determined by all these components together,
in fact one has \begin{equation}
\nor{a}_{\infty}=\prod_{p\in P}p^{v_{p}(a)}=\prod_{p\in P}\norp{a}^{-1}\qquad\forall a\in\gP.\label{eq-norm-reale-pidiche}\end{equation}

\section{Random walks on $\affP$ and $\mu$-boundaries}

\subsection{Random walk on $\affP$}

Consider now a probability measure $\mu$ on $\affP$ with a \emph{first
(logarithmic) moment} with respect to this geometry, that is \[
\int\sum_{p\in P}\nor{\ln\norp{a}}d\mu(a,b)<+\infty\quad\textrm{and }\quad\int\sum_{p\in\overline{P}}\ln^{+}\norp{b}d\mu(a,b)<+\infty\]
 and the associated random walk $R_{n}$ obtained as product of a
sequence $\left\{ g_{n}=(a_{n},b_{n})\right\} _{n}$of random elements
of $\affP$. A simple calculation shows that \[
R_{n}=(A_{n},Z_{n})=(a_{1}\cdots a_{n},\sum_{k=1}^{n}a_{1}\cdots a_{k-1}b_{k}).\]
The decisive parameters to control its behavior are the \emph{$p$-drifts}
\[
\phi_{p}=\int_{\affQ}\ln\norp{a}d\mu(a,b)=\EE(\ln\norp{a_{1}})\qquad\forall p\in\overline{P},\]
 that determine on which of the $p$-adic fields the action is contracting
in mean.

\subsection{$\mu$-boundaries}

It is known that whenever a $p$-drift is negative the corresponding
field $\QQ_{p}$ is a $\mu$-boundary; in fact one has the following 

\begin{lem}
Suppose that $\mu$ has a first moment. If $\phi_{p}<0$ for some
$p\in\overline{P}$ then the infinite sum \begin{equation}
Z_{\infty}^{p}=\sum_{k=1}^{\infty}a_{1}\cdots a_{k-1}b_{k}\label{eq-Z-in-Qp}\end{equation}
converges almost surely in $\Qp$ to the a random element and, for
all $x\in\Qp$, \[
R_{n}\cdot x=A_{n}x+Z_{n}\,\rightarrow\, Z_{\infty}^{p}\qquad\PP-\textrm{almost surely in }\Qp.\]

\end{lem}
\begin{proof}
For reader convenience, we give a sketch of the proof. For more details
see, for instance, \cite{El82} for the real case and \cite{CKW}
for the ultra-metric case.

Observe that by the Law of large numbers \[
\norp{a_{1}\cdots a_{n}}=\exp\left(\sum_{i=1}^{n}\ln\norp{a_{i}}\right)\]
converges to zero with exponential speed, roughly as $\exp(n\phi_{p})$.
On the other hand, since $\ln^{+}\norp{b_{1}}$ is integrable, $\left.\ln^{+}\norp{b_{n}}\right/n$
converges almost surely to zero. Thus the infinite sum (\ref{eq-Z-in-Qp})
converges, because its general term goes to zero exponentially. 
\end{proof}
Thus if $\nu_{p}$ is the law of $Z_{\infty}^{p}$ then for every
bounded function $\psi$ on $\QQ_{p}$\[
f(a,b)=\int_{\QQ_{p}}\psi(ax+b)d\nu_{p}(x)\]
 is a bounded harmonic function on $\affP$ and $(\Qp,\nu_{p})$ is
a quotient of the Poisson boundary. But it also follows almost immediately
that \[
R_{n}\cdot(x_{p})_{p}=\left(A_{n}x_{p}+Z_{n}\right)_{p}\,\rightarrow\,(Z_{\infty}^{p})_{p}\qquad\PP-\textrm{almost surely in }\prod_{p\in\overline{P}:\phi_{p}<0}\Qp\]
 for all $(x_{p})_{p}\in\prod_{p\in\overline{P}:\phi_{p}<0}\Qp.$ 

\begin{cor}
The Cartesian product of all the fields $\Qp$ whose drift is strictly
negative, $B=\prod_{p\in\overline{P}:\phi_{p}<0}\Qp$, endowed with
the joint law $\nu$ of $(Z_{\infty}^{p})_{p}$ is a $\mu$-boundary
of $\affP$.
\end{cor}
One can also interpret the $\mu$-boundary as an exit space of the
random walk. From this point of view the space $B$ together with
an extra point $\varpi$ provides the boundary of a topological compactification
of $\affP$ in which a sequence $\gamma_{n}=(\alpha_{n},\beta_{n})$
converges to a point $(x_{p})_{p}$ of $B$ if and only if \[
\alpha_{n}\rightarrow0\quad\textrm{and}\quad\beta_{n}\rightarrow x_{p}\quad\textrm{in }\Qp\,\textrm{for all }p\in\overline{P}\,\textrm{ such that }\phi_{p}<0.\]
Then the random walk, which is transient on the group $\affP$, converges
to a random element in $B$ with law $\nu$. 

Observe the the $p$-drifts are not all indipendent, but by (\ref{eq-norm-reale-pidiche})
they have to satisfy the following relation\[
\phi_{\infty}=-\sum_{p\in P}\phi_{p}\]
and thus they cannot all\- Hence, the diagonal embedding of $\ZP$
in $B$ is always dense, contrary to what one has in Lemma \ref{lem-appro-raggioI},
and the action of $\affP$ on $B$ has dense orbits. It follows that,
if the semi-group generated by the support $\mu$ is $\affP$, then
the support of the measure $\nu$ is the whole of $B$.

\subsection{A graphical interpretation}

Before going on, we would like to give a graphical interpretation
of the algebraic objects we are dealing with. Even if this is not
essential later on, this point of view was very useful in our preparatory
work and it can improve the geometrical understanding. 

First of all, as $\affP$ is a sub-group of the real affine transformation,
our group acts by isometries on the hyperbolic half-plane\[
\HH=\left\{ (x,y)|x\in\Rp,y\in\RR\right\} \]
by $(a,b)\cdot(x,y)=(\nor{a}x,ay+b)$.

\begin{figure}
\begin{center}\includegraphics[%
  width=0.70\textwidth]{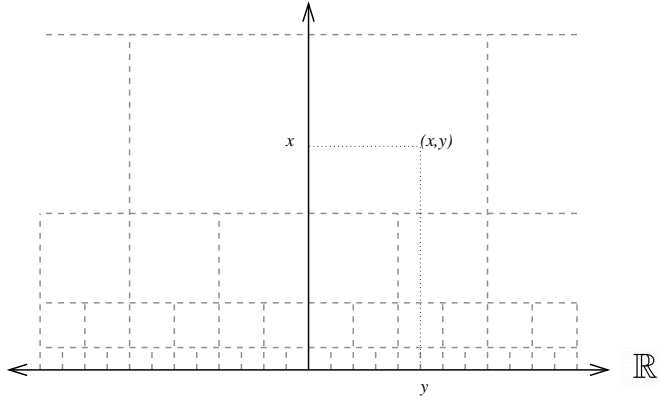}\end{center}

\caption{The hyperbolic half-plane}
\end{figure}

The $p$-adic equivalent of this structure is the oriented homogeneous
tree $\TT_{p}$ with degree $p+1$. This is the discrete graph without
loops and where each vertex has exactly $p+1$ neighbors, which can
be constructed as follows. Let $D^{p}(k,z)$ be the disc of radius
$p^{k}$ and center $z$ in $\QQ_{p}$. Since the valuation $v_{p}$
is integer valued and the $p$-adic norm has the ultra-metric property,
the set of all discs $V=\left\{ D^{p}(k,z)\right\} _{k\in\ZZ,z\in\Qp}$
is countable. Let us now define a graph $\TT_{p}$ whose vertices
are the discs and in which $D^{p}(k,z)$ is neighbor of $D^{p}(k+1,z)$
and $D^{p}(k-1,z)$. Since any disc of radius $p^{k}$ contains exactly
$p$ discs of radius $p^{k-1}$ and is contained in only one disc
of radius $p^{k+1}$, the graph $\TT_{p}$ is an homogeneous tree
of degree $p+1$. This tree is naturally devised in levels, indexed
by $\ZZ$, that consist of the discs of same radius. Furthermore,
any $z\in\Qp$ can be associated with the element of the boundary
of the tree, $\partial\TT_{p}$, that corresponds to the infinite
geodesic $\left\{ D^{p}(-k,z)\right\} _{k\in\NN}$, while we denote
by $\omega_{p}$ the end of the geodesic $\left\{ D^{p}(k,z))\right\} _{k\in\NN}$.
In this way one constructs an homeomorphism between $\Qp$ and the
boundary $\partial\TT_{p}-\left\{ \omega_{p}\right\} $ endowed with
the topology induced by the tree. As our affinities map a disc onto
a disc, one can define the action \[
(a,b)\cdot D^{p}(k,z)=D^{p}(k-v_{p}(a),az+b)\]
of $\affP$ on $\TT_{p}$, which is in fact an action by isometries
that fix the end $\omega_{p}$, similarly to the real case.

The simultaneous action of $\affP$ on the Cartesian product of $\HH$
and of all the $\TT_{p}$ for $p\in P$ has then discrete orbits by
Lemma \ref{lem-appro-raggioI} and corresponds to the embedding into
$\gP\times\prod_{p\in\overline{P}}\QQ_{p}$, when one neglects the
sign of the contracting-dilating coefficient. In fact fix as the origin
in $\HH\times\prod_{p\in P}\TT_{p}$ the point $o=\left((1,0),\left(D^{p}(0,0)\right)_{p\in P}\right)$.
Then the stabilizer of $o$ in $\affP$ is $\left\{ (1,0),(-1,0)\right\} $,
and, up to the sign, we can identify the group with the discrete orbit
\[
\affP o=\left\{ \left((\nor{a},b),\left(D^{p}(-v_{p}(a),b)\right)_{p\in P}\right)\in\HH\prod_{p\in P}\TT_{p}\left|a\in\gP,b\in\ZP\right.\right\} .\]
Observe that there is a relation between the level in the hyperbolic
plane and the heights on the $p$-adic trees that is given by the
relation between the $p$-norms (\ref{eq-norm-reale-pidiche}).

\begin{figure}
\begin{center}\includegraphics[%
  width=0.70\textwidth]{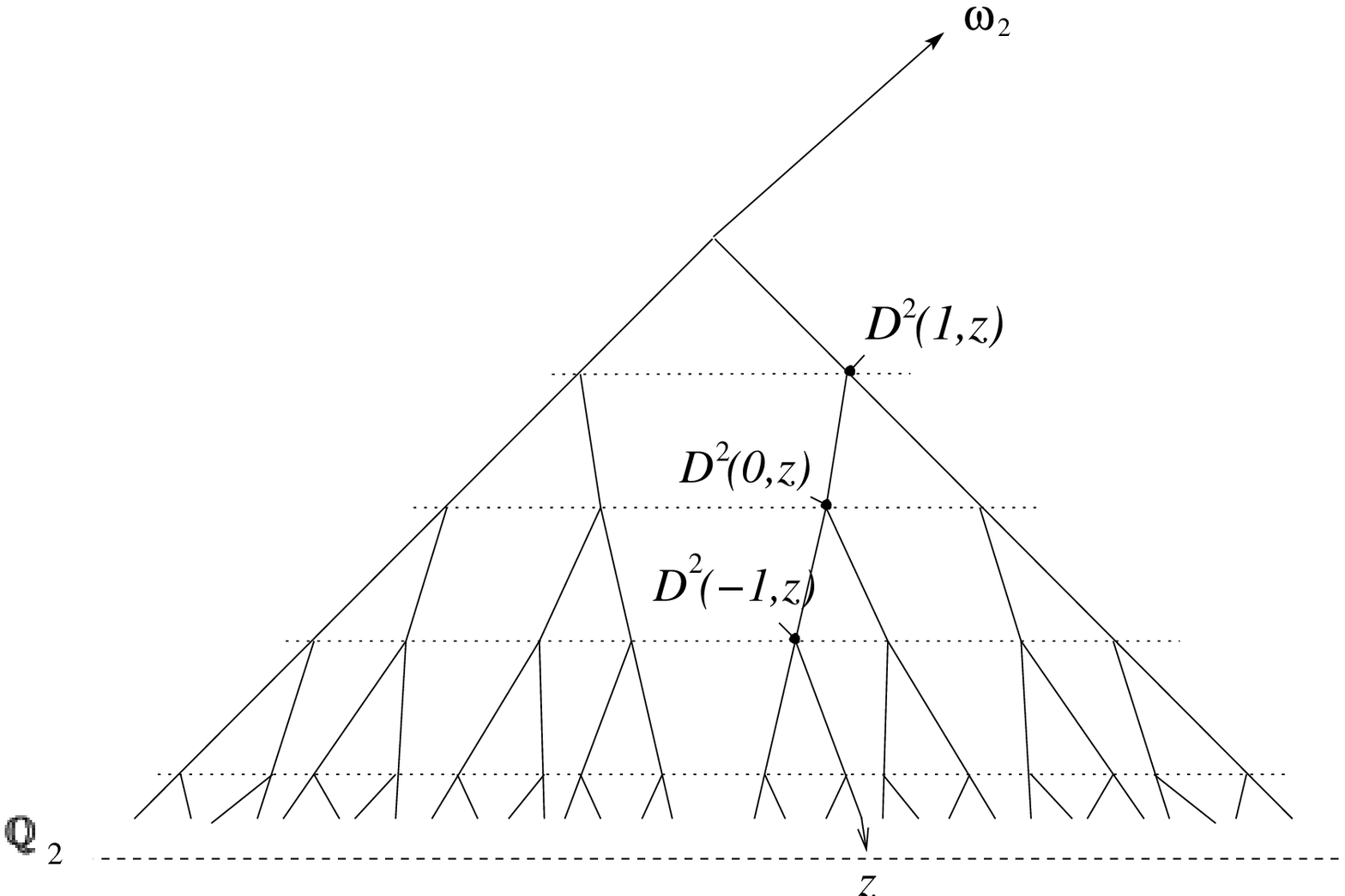}\end{center}

\caption{The tree of $\QQ_{2}$}
\end{figure}

In the case of the Baumslag-Solitar group this interpretation corresponds
to the one given by Farb and Mosher in \cite{FM99}. Here one has
just the half-plane and one tree and the link between the norms just
says that the level on $\HH$ is the opposite to the level on $\TT_{p}$.

\section{Poisson boundary}

\subsection{The strip criterion}

In the previous section we proved that the Cartesian products of the
$p$-adic fields on which the random walk has a contracting action
is a $\mu$-boundary and is a good candidate to be the maximal one,
that is the Poisson boundary. We will prove this fact using the the
geometrical strip criterion of Kaimanovich developed in \cite{Ka00}
(partially based on an idea of Ballmann and Ledrappier \cite{BL94}
and relying on the relationship between the boundary and the entropy
of the measure introduced by Avez \cite{Ave74} and developed by Vershik
and Kaimanovich \cite{VK79} and Derrienic \cite{Der80}). 

The starting point is to extend the random walk $R_{n}$ that is defined
usually only for positive times $n\in\NN$ also to the past using
a bilateral sequence of i.i.d. random variables $\left\{ g_{n}\right\} _{n\in\ZZ}$
on $G$ and the recursive equation \[
R_{n+1}=R_{n}g_{n+1}\qquad\textrm{and }\qquad R_{0}=e\]
to all $n\in\ZZ$, obtaining a bilateral path $\left\{ R_{n}\right\} _{n\in\ZZ}$.
Observe that for negative times $-m$, \[
R_{-m}=g_{0}\inv\cdots g_{-m+1}\inv=\check{R}_{m}\]
 is a random walk on $G$ whose law $\check{\mu}$ is the image of
$\mu$ under the inversion of the group. If $(B,\nu)$ and $(\check{B},\check{\nu})$
are a $\mu$-boundary and a $\check{\mu}$-boundary respectively,
one can imagine the trajectories of the bilateral process $\left\{ R_{n}\right\} _{n\in\ZZ}$
as paths in the group that join the exit point $\check{Z}_{\infty}$
of the random walk $\check{R}_{n}$ in $\check{B}$ to the exit point
$Z_{\infty}$ of $R_{n}$ in $B$. In order to prove that a $\mu$-boundary
is a Poisson boundary, one has to find some sort of strips in the
group connecting $\check{Z}_{\infty}$ to $Z_{\infty}$ that are visited
fairly often by the process $\left\{ R_{n}\right\} _{n\in\ZZ}$, without
being to large. 

One way to build sets that verify the first of these conditions is
to give a map $S$ that associates with every couple of exit points
$(\check{z},z)\in\check{B}\times B$ a non-empty subset of the group
$G$ that is equinvariant for the action of the group:\[
gS(\check{z},z)=S(g\check{z},gz)\qquad\forall g\in G.\]
 In Theorem 6.4 of \cite{Ka00}, Kaimanovich shows that, when $\check{z}$
and $z$ are distributed according to the exit measures, the probability
for $R_{n}$ being in $S(\check{z},z)$ is constant for all $n\in\ZZ$
and that one has the following criterion:

\begin{prop}
\label{prop-striscie-doppie} Suppose that the measure $\mu$ has
finite first moment. If there exists a sequence of sets $C_{n}$ such
that:\[
\pr{R_{n}\in C_{n}}>\varepsilon\qquad\forall n\in\NN\]
for some $\varepsilon>0$, and $\check{\nu}(d\check{z})\times\nu(dz)$-almost
surely \[
\lim_{n\toinf}\frac{1}{n}\ln\card{S(\check{z},z)g\cap C_{n}}=0\qquad\forall g\in G,\]
 then $B$ is the Poisson boundary.
\end{prop}

\subsection{Strips in $\affP$}

Let now go back to the group $G=\affP$. In the previous section we
have shown that \[
B=\prod_{p\in\overline{P}:\phi_{p}<0}\Qp\]
is a $\mu$-boundary. Since the $p$-drifts of the reflected measure
$\check{\mu}$ are\[
\check{\phi}_{p}=\esp{\ln\norp{a_{1}^{-1}}}=-\phi_{p},\]
 the space \[
\check{B}=\prod_{p\in\overline{P}:\check{\phi}_{p}<0}\Qp=\prod_{p\in\overline{P}:\phi_{p}>0}\Qp\]
 is the right candidate for being the maximal $\check{\mu}$-boundary.

For every point $z_{p}$ of $\Qp$ we define the {}``cone'' in $G$
of vertex $z_{p}$ as \[
C^{p}(z_{p})=\left\{ (a,b)\in G:\norp{z_{p}-b}\leq\norp{a}\right\} =\left\{ g\in G:\norp{g^{-1}\cdot z}\leq1\right\} .\]
 When $p=\infty$ and one embeds the group $G$ in the hyperbolic
half-plane, this is, the set of points that are between the two lines
that start at the point of coordinate $(0,z_{\infty})$ and have a
slope of 45 degrees. For $p\neq\infty$, if one projects the cone
on the corresponding tree, one obtains the geodesic that connects
$z_{p}$ to the end $\omega_{p}$.

We can now define the strip from $\check{\mathbf{z}}=(z_{p})_{p:\phi_{p}>0}\in\check{B}$
to $\mathbf{z}=(z_{p})_{p:\phi_{p}<0}\in B$ as\[
S(\check{\mathbf{z}},\mathbf{z})=\bigcap_{p\in\overline{P}:\phi_{p}\neq0}C^{p}(z_{p})=\left\{ (a,b)\in G:\norp{z_{p}-b}\leq\norp{a}\,\,\forall p\in\overline{P}:\phi_{p}\neq0\right\} .\]
These strips are equivariant as the cones $C^{p}(z_{p})$.%
\begin{figure}
\begin{center}\includegraphics[%
  width=1.0\textwidth]{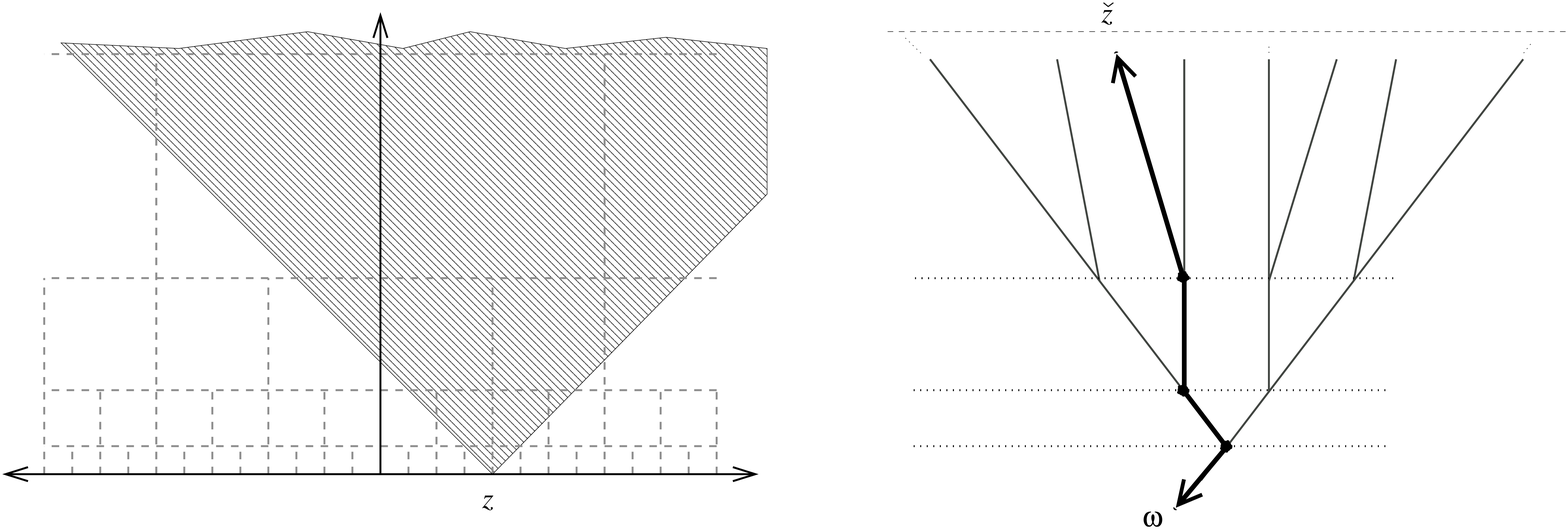}\end{center}

\caption{Strip in $\affO{\left\{ 2\right\} }$ or $BS(1,2)$}
\end{figure}
 In the case of the Baumslag-Solitar group $BS(1,p)$, when the real
drift is negative, the boundary $B$ is the real line while $\check{B}$
is $\QQ_{p}$. Using its discrete embedding in $\HH\times\TT_{p}$,
a strip is the set consisting of the pairs of a point of the cone
in the plane and a vertex of the geodesic in the tree that are at
opposite levels.

When all the $p$-drifts are null, this still applies even thought
in a degenerate way: in fact the boundaries $B$ and $\check{B}$
are trivial and collapse to a single point and the strips become the
whole group $G$. 

Observe that if for all $p\in\overline{P}$ the drift $\phi_{p}$
is not zero, in the definitionof a strip $S(\check{\mathbf{z}},\mathbf{z})$
all the significant prime $p$ are involved. In this case the intersection
of a strip with the translations group, that is \[
\pi(\check{\mathbf{z}},\mathbf{z})=S(\check{\mathbf{z}},\mathbf{z})\cap Tra(\ZP)=\left\{ (1,b):b\in\ZP\textrm{ and }\norp{z_{p}-b}\leq1\,\,\forall p\in P\right\} ,\]
 has alway 2 or 3 elements (by Lemma\ref{lem-appro-raggioI}). On
the other hand \[
S(\check{\mathbf{z}},\mathbf{z})=\bigcup_{a\in\gP}(a,1)\pi(a\inv\check{\mathbf{z}},a\inv\mathbf{z}).\]
 As the group $\gP$ has polynomial growth and the projection of the
random walk on $\gP$ grows linearly, we can directly apply the strip
criterion of Proposition \ref{prop-striscie-doppie} with \[
C_{n}=\left\{ (a,b)\in G:\max_{p\in P}\left|\ln\norp{a}\right|\leq Kn\right\} \]
for a suitable $K>0$.

\subsection{Proof of Theorem \ref{thm-poisson-boundary}}

When some of the drifts are centered, in order to conclude one needs
to show that the growth of the random walk in the centered directions
is sub-exponential, using the following lemma

\begin{lem}
Let $Z_{n}=\sum_{1\leq k\leq n}a_{1}\cdots a_{k-1}b_{k}$. For every
$p\in\overline{P}$ such that $\phi_{p}=0$, one has \[
\lim\frac{1}{n}\ln^{+}\norp{Z_{n}}=0\qquad\PP-\textrm{almost surely}.\]

\end{lem}
\begin{proof}
Let $\norp{a_{1}\cdots a_{n}}=\e^{S_{n}}$. It follows straightly
from the strong Law of large numbers that $S_{n}/n$ converge almost
surely to zero. On the other hand \[
\left|Z_{n}\right|_{p}=\norp{\sum_{1\leq k\leq n}a_{1}\cdots a_{k-1}b_{k}}\leq n\left(\max_{1\leq k\leq n}\norp{a_{1}\cdots a_{k-1}b_{k}}\right).\]
As the random variables $\ln^{+}\norp{b_{n}}$ are integrable, $(\ln^{+}\norp{b_{n}})/n$
converge to zero and we can deduce that\[
\limsup_{n\toinf}\frac{\ln\left|Z_{n}\right|_{p}}{n}\leq\limsup_{n\toinf}\frac{\ln n}{n}+\frac{\max_{1\leq k\leq n}\left\{ S_{k}+\ln^{+}\norp{b_{k}}\right\} }{n}\leq0\]

\end{proof}
We have introduced all the tools we need and we can finally give the
proof of our main result 

\begin{proof}
of Theorem \ref{thm-poisson-boundary}. Observe that by the previous
lemma there exists a sequence $\varepsilon_{n}$ that converges to
zero such that \[
\pr{\max_{p\in\overline{P}:\phi_{p}=0}\ln^{+}\norp{Z_{n}}\leq\varepsilon_{n}n}\,\rightarrow\,1\]
and that by the Law of large numbers \[
\lim_{n\toinf}\frac{1}{n}\max_{p\in P}\left|\ln\norp{A_{n}}\right|=\max_{p\in P}\left|\phi_{p}\right|<\infty\]
$\PP$-almost surely. For any $K>\max_{p\in P}\left|\phi_{p}\right|$,
let \[
C_{n}=\left\{ (a,b)\in G:\max_{p\in P}\left|\ln\norp{a}\right|\leq Kn\textrm{ and }\max_{p\in\overline{P}:\phi_{p}=0}\ln^{+}\norp{b}\leq\varepsilon_{n}n\right\} .\]
 Then $\pr{R_{n}\in C_{n}}$ converges to 1. 

Let now \[
Q_{n}=\left\{ (a,b)\in\affP:\max_{p\in P}\left|\ln\norp{a}\right|\leq Kn\textrm{ and }\max_{p\in\overline{P}}\ln^{+}\norp{b}\leq\varepsilon_{n}n\right\} ,\]
which is the analogue of $C_{n}$ but where one considers all the
$p$-norms. Consider the section\[
H=\left\{ (1,b)\in\affP:\norp{b}\leq1\,\,\forall p\in\overline{P}:\phi_{p}=0\right\} .\]
One can easily check the following relation\[
C_{n}\subseteq Q_{n}H.\]
By Lemma \ref{lem-appro-raggioI}, one knows that the cardinality
of the set \[
S(\check{z},z)\cap H=\left\{ (1,b)\in\affP:\norp{b-z_{p}}\leq1\,\forall p:\phi_{p}\neq0\textrm{ and }\norp{b}\leq1\,\forall p:\phi_{p}=0\right\} \]
is less then 3 for all $(\check{z},z)\in\check{B}\times B$. Thus
we can estimate the cardinality of the strip along the sets $C_{n}$
using the equinvariance \[
\textrm{card}(S(\check{z},z)\cap C_{n})\leq\textrm{card}Q_{n}\cdot\sup_{(\check{z},z)\in\check{B}\times B}\textrm{card}(S(\check{z},z)\cap H)\leq3\textrm{card}Q_{n}.\]
 Finally to control the cardinality of $Q_{n}$, observe that for
any $m>1$ \[
\card{b=h\prod_{p\in P}p^{h_{p}}\in\ZP:\norp{b}\leq m\,\forall p\in\overline{P}}\leq c_{1}\e^{2m}\]
 and \[
\card{a=\pm\prod_{p\in P}p^{k_{p}}\in\gP:\norp{a}\leq m\,\forall p\in\overline{P}}\leq m^{c_{2}}\]
for some constants $c_{1}$ and $c_{2}$ depending on $P$. Thus one
has\[
\frac{1}{n}\ln\textrm{card}(S(\check{z},z)\cap C_{n})\leq\frac{\ln3+c_{2}\ln(Kn)+\ln c_{1}+2\varepsilon_{n}n}{n}\,\rightarrow\,0.\]
 To conclude just observe that for any $g\in\affP$ there exits a
constant $h$ such that \[
\textrm{card}(S(\check{z},z)g\cap C_{n})=\textrm{card}(S(\check{z},z)\cap C_{n}g^{-1})\leq\textrm{card}(S(\check{z},z)\cap C_{n+h}).\]
 Thus all hypothesis of Proposition \ref{prop-striscie-doppie} are
verified.
\end{proof}
\bibliographystyle{alpha}
\bibliography{sara}

\lyxrightaddress{Sara BROFFERIO\\
 Institut fuer Mathematik C\\
 Technische Universitaet Graz\\
 Steyergasse 30\\
 A-8010 Graz\\
 brofferio@finanz.math.tu-graz.ac.at}
\end{document}